\documentclass[12pt]{article} 
\usepackage{amsfonts,amssymb}
\usepackage{exscale}
\usepackage[notcite,notref]{}

\usepackage{amsgen,amsmath,amstext,amsbsy,amsopn,amsfonts,amssymb,pifont, epsfig}
\usepackage{graphicx}

\let\BFseries\bfseries\def\bfseries{\BFseries\mathversion{bold}} 

\newtheorem{thm}{Theorem}

\newtheorem{cor}{Corollary}

\newenvironment{proof}[1][] {\smallskip \noindent {\bf Proof#1.} }{\hspace*{\fill}$\square$\medskip\par}

\newcommand{\be}{\begin{equation}}
\newcommand{\ee}{\end{equation}}
\newcommand{\bea}{\begin{eqnarray}}
\newcommand{\bq}{\begin{eqnarray}}
\newcommand{\eea}{\end{eqnarray}}
\newcommand{\eq}{\end{eqnarray}}
\newcommand{\beaa}{\begin{eqnarray*}}
\newcommand{\eeaa}{\end{eqnarray*}}
\newcommand{\bei}{\begin{itemize}}
\newcommand{\eei}{\end{itemize}}
\newcommand{\bee}{\begin{enumerate}}
\newcommand{\eee}{\end{enumerate}}
\newcommand{\bi}{\begin{itemize}}
\newcommand{\ei}{\end{itemize}}
\newcommand{\beq}{\begin{eqnarray*}}
\newcommand{\eeq}{\end{eqnarray*}}
\newcommand{\beqn}{\begin{eqnarray}}
\newcommand{\eeqn}{\end{eqnarray}}

\newcommand{\ignore}[1]{}{}

\newcommand{\e}{\varepsilon}
\newcommand{\eps}{\varepsilon}

\def\<{\langle}
\def\>{\rangle}

\def\norm#1{\left\|#1\right\|}             

\def\P{{\mathbb P}}
\def\pr#1{\P\left(#1\right)}

\def\E{{\mathbb E}\,}

\def\R{{\mathbb R}}

\def\N{\mathbb N}

\newcommand{\dd}{\, {\rm d}}

\begin{document}
\begin{titlepage}


\title{Small deviations for a family of smooth Gaussian processes}
\author{Frank Aurzada\thanks{Institut f\"{u}r Mathematik, Technische Universit\"{a}t Berlin, Germany, \texttt{aurzada@math.tu-berlin.de}, supported by the DFG Emmy Noether program.} \and
Fuchang Gao\thanks{
Department of Mathematics, University of Idaho, \texttt{fuchang@uidaho.edu}, supported in part by NSF grant DMS-0806211.}
\and
Thomas K\"{u}hn\thanks{Mathematisches Institut, Universit\"{a}t Leipzig, Germany, \texttt{kuehn@math.uni-leipzig.de}, supported in part by Ministerio de Ciencia e Innovaci\'{o}n (Spain), grant MTM2010-15814.}
\and
Wenbo V. Li\thanks{
Department of Mathematical Sciences,
University of Delaware,
\texttt{wli@math.udel.edu}, supported in part by NSF grant DMS--0805929, NSFC-6398100, CAS-2008DP173182.}
\and
 Qi-Man Shao\thanks{Department of Mathematics, Hong Kong University of Science and Technology, Hong Kong, \texttt{maqmshao@ust.hk},
supported in part by Hong Kong RGC CERG 602608.}
}
\date{\today}
\maketitle
\thispagestyle{empty}

\begin{abstract}
We study the small deviation probabilities of a family of very smooth self-similar Gaussian processes.
The canonical process from the family has the same scaling property as standard Brownian motion
and plays an important role in the study of zeros of random polynomials.

Our estimates are based on the entropy method, discovered in Kuelbs and Li (1992) and developed further in Li and Linde (1999), Gao (2004), and Aurzada et al.\ (2009). While there are several ways to obtain the result w.r.t.\ the $L_2$ norm, the main contribution of this paper concerns the result w.r.t.\ the supremum norm. In this connection, we develop a tool that allows to translate upper estimates for the entropy of an operator mapping into $L_2[0,1]$ by those of the operator mapping into $C[0,1]$, if the image of the operator is in fact a H\"older space.

The results are further applied to the entropy of function classes, generalizing results of Gao et al.\ (2010).
\end{abstract}



\end{titlepage}

\setcounter{page}{1}

\section{Introduction and main results}
\subsection{Introduction}
The small deviation problem for a sto\-chas\-tic process $X=(X(t))_{t\geq 0}$ -- also called small ball or small value problem -- consists in determining the probability
$$
-\log \pr{ ||X|| \leq \eps}, \qquad\text{as $\eps\to 0$,}
$$
where $||.||$ is for example the norm in some $L_p[0,1]$ or in $C[0,1]$. Small deviation probabilities play a fundamental role in many problems in probability and analysis, see the lecture notes \cite{Li10} for details. This is why there has been a lot of interest in small deviation problems in recent years, cf.\ the survey \cite{LS} and the literature compilation \cite{sdbib}. There are many connections to other questions such as the law of the iterated logarithm of Chung type, strong limit laws in statistics, metric entropy properties of linear operators, quantization, and several other approximation quantities for stochastic processes. For Gaussian processes, a reasonable amount of theory has been developed up to date, see e.g.\ \cite{LS}.

In this paper, we study a family of very smooth self-similar Gaussian processes and their respective small deviation probabilities. This family of processes  does not seem to have drawn enough attention in the probability community. Yet, it is very important and appears naturally in physics \cite{Majumdar}.

Let us give a short motivation for studying this class of processes. Let $\phi({\bf x}, t)$ be a scalar field in a $d$-dimensional space that evolves according to the
deterministic diffusion equation 
 $$\frac{\partial}{\partial t}\phi({\bf x}, t)=  \nabla^2 \phi({\bf x}, t)$$
with the initial condition $\phi({\bf x}, 0)=\psi({\bf x})$. The only randomness is the initial condition which could be a mean-zero Gaussian random field.
For example, $\phi({\bf x}, t)$ could represent the density fluctuation of a diffusing gas. In such a case, it is reasonable to assume that $\psi({\bf x})$ is of zero mean and delta covariance $\E \psi({\bf x})\psi({\bf y})=\delta({\bf x-y})$.
For a system of linear size $L$, the solution of the diffusion equation in the bulk of the system is
$$\phi({\bf x}, t)=(4\pi t)^{-d/2}\int_{|{\bf y}|\le L}\exp(-\|{\bf x-y}\|_2^2/4t)\psi({\bf y}) \dd {\bf y}.$$
Because of the linearity of the integral, $\phi({\bf x}, t)$, $t\ge 0$, is a Gaussian process. It is customary to study the normalized process $X(t)=\phi({\bf x}, t)/(\E \phi({\bf x}, t)^2)^{1/2}$. When $L=\infty$, it is straightforward to calculate that the Gaussian process $X(t)$ has covariance structure 
$$\E X(t)X(s)=\left(\frac{2ts}{(t+s)^2}\right)^{d/4}.$$
In modeling non-Fickian diffusion, it is common that the Green function $$G({\bf x},t)=(4\pi t)^{-d/2}\exp(-\|{\bf x}\|_2^2/4t)$$ above needs to be replaced by a different Green function with a different scaling property, say $G_{p,q}({\bf x},t)=Ct^{-p}K({\bf x}/t^q)$.

This naturally leads us to consider a centered Gaussian processes $X_{\alpha, \beta}(t)$, $ t > 0$,
defined by the covariance function
\be\label{def:K}
K(t,s)=\E X_{\alpha, \beta}(t) X_{\alpha, \beta}(s) = \frac{2^{2\beta+1} (ts)^{\alpha}}{(t+s)^{2\beta+1}}, \quad t,s>0,
\ee
for $\alpha>0$ and $\beta >-1/2$. Note that for $\alpha > \beta+1/2$, we can define $X_{\alpha, \beta}(0)=0$. 
It is also easy to check that $X_{\alpha, \beta}$ is an $(\alpha-\beta -1/2)$-self-similar process,
i.e.\ $(X_{\alpha, \beta}(ct))$ has the same law as $(c^{\alpha-\beta -1/2} X_{\alpha, \beta}(t))$ for any $c>0$.
In particular, $X_{\alpha, \beta}$ has the same scaling property as Brownian motion for $\alpha-\beta=1$.
A useful stochastic integral representation is
\be\label{def:x}
X_{\alpha, \beta}(t) = \sqrt{ \frac{2^{2\beta+1}}{\Gamma(2\beta+1)} } \, t^{\alpha} \int_0^\infty x^{\beta} e^{-x t} \dd B(x),\qquad t > 0,
\ee
where $B$ is a standard Brownian motion.

If $\beta=0$, it is easy to see, using integration by parts, that $X_{\alpha,0}$ has the same law as the process
$$
{\tilde X}_{\alpha, 0}(t) :=  \sqrt{2}\, t^{1+\alpha} \int_0^\infty e^{-x t} B(x) \, \dd x,\qquad t > 0.
$$

We also introduce, for $\alpha=1$, $\beta=0$, the canonical process
$$
X(t):=X_{1,0}(t)= \sqrt{ 2}\, t \int_0^\infty  e^{-x t} \dd B(x),\qquad t \ge 0\,.
$$
This process has the same scaling property as standard Brownian motion
and plays an important role in the study of zeros of random polynomials, cf.\ \cite{LS04,LS05,ad09}, which gives another reason for studying this process. Yet a different motivation is the applicability of small deviation estimates of stationary processes in statistics, cf.\ e.g.\ \cite{vaartzanten}; we come back to this point in Section~\ref{sec:lamperti} below.

This paper is structured as follows. In Section~\ref{sec:mainres}, we present our main results  concerning the small deviation rate w.r.t.\ $L_2[0,1]$ and supremum norm. We compare the present findings to other known facts for smooth Gaussian processes in Section~\ref{sec:lamperti}. The main tool for the proofs is the so-called entropy method, which is recalled in Section~\ref{sec:entropymethod}. We also present the functional analytic counterparts of our main theorems in that section.

There are different ways to obtain the small deviation rate for the $L_2$-norm for our class of processes. Therefore, the main contribution of this paper is the transfer from $L_2$ to supremum norm estimates. One interesting tool in this connection is presented in Section~\ref{sec:holder}, and we believe that it may be of independent interest.

The proofs of all statements are given in Section~\ref{sec:proofs}. Finally, our main theorems entail some new results for the entropy of function classes. This connection and the respective corollaries are presented in Section~\ref{sec:functionS}.

We remark that further properties of the class of processes $(X_{\alpha,\beta})$ and alterative proofs are contained in an earlier version of this paper available from http://arxiv.org/abs/1009.5580.

Let us fix some notation. We write $f \preceq g$ or $g \succeq f$ if $\limsup f/g < \infty$,
and the asymptotic equivalence $f \asymp g$ means that we have both $f \preceq g$
and $g \preceq f$. Moreover, we write $f \lesssim g$ or $g \gtrsim f$, if $\limsup f/g \le  1$.
Finally, the strong equivalence $f \sim g$ means that $\lim f/g = 1$.

\medskip
{\bf Acknowledgement:}  We are grateful to the American Institute of Mathematics (AIM) for supporting the workshop ``Small ball inequalities in analysis, probability, and irregularities of distribution'' (December 2008), where the work on this paper was started. We would like to thank two anonymous referees and an associate editor for their suggestions, that helped to improve the exposition of the paper.

\subsection{Main results}\label{sec:mainres}
We can clarify the small deviation order for the whole class of smooth processes introduced above for the $L_2$-norm and the sup-norm.

For the $L_2$-norm, there is the following result.

\begin{thm} \label{thm:l2smallballs}
Let $\alpha>\beta>-1/2$. Let $X_{\alpha, \beta}$ be the process defined by (\ref{def:K}).
Then $$-\log \pr{ \int_0^1 |X_{\alpha, \beta} (t)|^2 \, \dd t \leq \eps^2} \sim \kappa_{\alpha,\beta} |\log \eps|^3,$$
where the constant is given by
\begin{equation} \label{eqn:sdconst2}
\kappa_{\alpha,\beta} := \frac{1}{3 (\alpha-\beta)\pi^2}.
\end{equation}
\end{thm}

We remark that Theorem~\ref{thm:l2smallballs} follows from Proposition~4.3 in \cite{kn10}, which we got to know only after the submission of the present paper. We will outline their method of proof and the relation to our approach below.

For the sup-norm, we obtain the following theorem under optimal assumptions on the parameters $\alpha$ and $\beta$. However, the result is less precise with respect to the small deviation constant compared to the $L_2$-norm result.

\begin{thm} \label{thm:linftysmallballs}
Let $\alpha>\beta+1/2>0$. Let $X_{\alpha, \beta}$ be the process defined by (\ref{def:K}) with $X_{\alpha, \beta}(0)=0$.
Then, with some constant $\tilde{\kappa}_{\alpha,\beta}>0$, we have
\begin{equation} \label{eqn:supnormcase}
\tilde\kappa_{\alpha,\beta} |\log \eps|^3 \lesssim -\log \pr{ \sup_{t\in[0,1]} |X_{\alpha, \beta} (t)| \leq \eps} \lesssim \kappa_{\alpha-1/2,\beta} |\log \eps|^3,
\end{equation}
where $\kappa_{\alpha,\beta}$ is defined in (\ref{eqn:sdconst2}) and $\tilde\kappa_{\alpha,\beta}>0$. Further, $\tilde{\kappa}_{\alpha,\beta}\to\infty$ when $\alpha-\beta\to 1/2$.
\end{thm}

The proofs of these two theorems use the connection between small deviations of Gaussian processes and the entropy numbers of a linear operator generating the process,
cf.\ \cite{kl}, \cite{ll}, \cite{ailz}. In fact, due to Corollaries~2.2 and~2.4 in \cite{ailz},
Theorem~\ref{thm:l2smallballs} and Theorem~\ref{thm:linftysmallballs} are equivalent to Theorem~\ref{entropy2} and Theorem~\ref{entropyinfty}, respectively, given in Section~\ref{sec:entropymethod}. Other interesting small deviation estimates for smooth Gaussian processes can be found in \cite{ailz}, \cite{ku}, and \cite{kn10}, and we compare them to the present results in the next section.

\subsection{Stationary version of our process} \label{sec:lamperti}
Let us shortly comment on a relation to \cite{ailz} and \cite{kn10}, where, among other things, small deviation probabilities of stationary Gaussian processes are considered.

Our process $X_{\alpha,\beta}$ is $(\alpha-\beta-1/2)$-self-similar. We can thus consider the Lamperti transformed stationary Gaussian process: $Y_{\beta}(t):=e^{-(\alpha-\beta-1/2) t} X_{\alpha,\beta}(e^t)$, $t\in \R$, which has the correlation function
$$
k(t)=\frac{2^{2\beta+1}}{(e^{t/2}+e^{-t/2})^{2\beta+1}} = \frac{1}{\cosh(t/2)^{2\beta+1}},\qquad t\geq 0.
$$
One can find that for the spectral measure we have (for all $\beta>-1/2$)
\begin{equation}\label{eqn:spectralmeasure}
 k(t) =: \int_{-\infty}^\infty e^{i t x} e^{-G(x)} \dd x,\qquad \text{with}\qquad G(x)\sim \pi x,\quad x\to\infty.
\end{equation}
Therefore, by Proposition~4.2(3) in \cite{kn10}, for some numeric constant $c>0$,
\begin{equation}\label{eqn:kn}
-\log \pr{ \int_0^1 |Y_\beta(t)|^2 \dd t \leq \eps^2} \sim c\, |\log \eps|^2,
\end{equation}
as for the process in \cite{ailz} with $G(x)=|x|$. In comparison, note that our main theorems concern the infinite time horizon for the stationary process: Namely, for the $L_2$ norm, Theorem~\ref{thm:l2smallballs} reads in terms of $Y_\beta$:
$$
-\log \pr{ \int_{-\infty}^0 |Y_\beta(t)|^2 e^{2 (\alpha-\beta) t} \dd t \leq \eps^2} \sim \kappa_{\alpha,\beta}\, |\log \eps|^3,
$$
which makes sense since $\alpha-\beta>0$.

In the same spirit, we can compare the results w.r.t.\ the sup-norm. Using (\ref{eqn:spectralmeasure}) with Lemma~2.3 in \cite{vaartzanten} and (\ref{eqn:kn}) one obtains for the finite time horizon
$$
-\log \pr{ \sup_{0\leq t\leq 1} |Y_\beta(t)| \leq \eps} \asymp |\log \eps|^2,
$$
while our main result for the sup-norm, Theorem~\ref{thm:linftysmallballs}, is again for the infinite time horizon and reads
$$
-\log \pr{ \sup_{-\infty<t\leq 0} |Y_\beta(t) e^{(\alpha-\beta-1/2)t}| \leq \eps} \asymp |\log \eps|^3,
$$
which makes sense since $\alpha-\beta-1/2>0$.

\subsection{The entropy method} \label{sec:entropymethod}
In this section, we recall the entropy method for obtaining small deviation probabilities and formulate the respective counterparts for Theorems~\ref{thm:l2smallballs} and~\ref{thm:linftysmallballs}.

In order to state and use the connection to the entropy numbers, let us first define the entropy numbers. 
For a linear operator $u : E \to F$ between Banach spaces $E$ and $F$ and $n\in\mathbb{N}$, the entropy numbers are defined as follows:
$$
e_n(u : E \to F):= \inf \left\{ \eps>0 \,:\, \exists f_1, \ldots, f_{2^{n-1}} \in F \enspace{\rm s.t.}\enspace
u(B_E)\subseteq\bigcup_{k=1}^{2^{n-1}}(f_k+B_F)\right\}\,,
$$
where $B_E$ and $B_F$ denote the closed unit balls in $E$ and $F$, respectively. For elementary properties and further information see e.g.\ \cite{Carl}. Since $u$ is compact if and only if $\lim_{n\to\infty} e_n(u)=0$, the decay rate of the entropy numbers is a measure for  the ``degree of compactness" of $u$.

It turns out that there is a close relation between the small deviation problem for a Gaussian process $X$ attaining values in $E$ and the entropy numbers of a compact operator $u : L_2(S) \to E$ related to $X$ through its characteristic function:
\begin{equation}
 \label{eqn:acop}
 \E e^{i \langle X, h\rangle } = \exp\left( - \frac{1}{2} \, \norm{u'(h)}_{L_2(S)}^2\right),\qquad h\in E',
\end{equation}
where $u' : E' \to L_2(S)$ is the dual operator and $(S,\mathcal{S},\lambda)$ is some measure space. Here, $(E,\norm{.})$ is some Banach space.

It can be checked easily that, up to an unimportant multiplicative constant, the process $X_{\alpha,\beta}$ defined in (\ref{def:K}) is related -- via (\ref{eqn:acop}) -- to the operator
\be\label{uft}
(u f)(t) = t^\alpha \int_0^\infty x^\beta e^{- xt} f(x) \, \dd x,\qquad f\in L_2[0,\infty).
\ee

Note that the process defined in (\ref{def:K}) can be considered with values in the Banach spaces $E=L_2[0,1]$ for $\alpha>\beta>-1/2$. 
In this $L_2$-setting, we obtain the precise behavior of the entropy numbers of the operator related to our process on the exponential scale.

\begin{thm}\label{entropy2}
Let $\alpha>\beta>-1/2$. Let $u : L_2[0,\infty) \to L_2[0,1]$ be the operator given by (\ref{uft}).
Then $$-\log e_n(u)\sim d_{\alpha,\beta} \,n^{1/3},$$ where
\begin{equation} \label{eqn:l2const}
d_{\alpha,\beta}:=( 3(\alpha-\beta)\pi^2 \log 2)^{1/3}.
\end{equation}
\end{thm}

In the sup-norm case we get the following result under the optimal assumption on the parameters $\alpha-\beta>1/2$, which is exactly when the process defined in (\ref{def:K}) is almost surely in $E=C[0,1]$.
\begin{thm}\label{entropyinfty}
Let $\alpha>\beta+1/2>0$.
Let $u : L_2[0,\infty) \to C[0,1]$ be the operator given by (\ref{uft})
Then $$\tilde{d}_{\alpha,\beta} \,n^{1/3} \lesssim -\log e_n(u)\lesssim d_{\alpha-1/2,\beta} \,n^{1/3},$$ where $d_{\alpha,\beta}$ is defined in (\ref{eqn:l2const}) and 
$$
\tilde{d}_{\alpha,\beta}:=\frac{\min(\alpha-\beta-1/2\,,\,1) }{1/2+\min(\alpha-\beta-1/2\,,\,1)} \cdot d_{\alpha,\beta}.
$$
\end{thm}

Note, in particular, that if $\alpha-\beta\to 1/2$, then
$$
\tilde{d}_{\alpha,\beta} \to 0 \qquad\text{and}\qquad d_{\alpha-1/2,\beta}\to 0.
$$

\subsection{A lemma relating entropy in $L_2$ and $L_\infty$ via H{\"{o}}lder continuity} \label{sec:holder}
In this section, we present an interpolation result which will be an essential tool in the proof of Theorem~\ref{entropyinfty}, and which might be of independent interest. 

Roughly speaking, it provides a technique for deriving upper entropy estimates for an operator $u:E\to C[0,1]$, where $E$ is a Banach space, from entropy estimates of $u: E\to L_2[0,1]$, 
i.e.\ the same operator, but considered as operator into a larger target space. 
The additional information we need for this argument is that $u$ should map $E$ even into a smaller space than $C[0,1]$, namely into a H\"older space $C_\lambda[0,1]$ for some $0<\lambda \le 1$. 
This space consists of all functions $f\in C[0,1]$ such
$$
\Vert f\Vert_{C_\lambda[0,1]}:=\sup_{0\le s<t\le 1}\frac{|f(t)-f(s)|}{|t-s|^\lambda} + \sup_{0\le t\le 1}|f(t)| <\infty\,.
$$
Moreover, $C_\lambda[0,1]$ is a Banach space under the norm $\Vert . \Vert_{C_\lambda[0,1]}$.

\begin{thm}\label{general}
Let $u$ be an operator from a Banach space $E$ into a H\"older space $C_\lambda[0,1]$ for some $0< \lambda \le 1$.
Then we have for all $k\in \N$
$$
e_k(u:E \to C[0,1])\le 2\,\Vert u:E \to C_\lambda [0,1]\,\Vert^{\frac{1}{2\lambda +1}}\,e_k(u:E \to L_2[0,1])^{\frac{2\lambda}{2\lambda + 1}}\,.
$$
\end{thm}

\section{Proofs of the theorems} \label{sec:proofs}
\subsection{Proof of Theorems~\ref{thm:l2smallballs} and~\ref{entropy2}}
There are many ways to obtain the small deviation estimates w.r.t.\ the $L_2$ norm. Even though Theorem~\ref{thm:l2smallballs} and Theorem~\ref{entropy2} are equivalent due to Corollaries~2.2 and~2.4 in \cite{ailz}, we will give proofs for both of the theorems: On the one hand, the proof of Theorem~\ref{thm:l2smallballs} is simpler, on the other hand, for consistency with the $C[0,1]$ case it seems reasonable to have a purely analytic proof of Theorem~\ref{entropy2}. Additionally, we will explicitly use Theorem~\ref{entropy2} during the proof of Theorem~\ref{entropyinfty}, both for upper and lower bound.

We remark that the proof of Theorem~\ref{thm:l2smallballs} is as outlined in Proposition~4.3 in \cite{kn10}, which we got to know only after the submission of the present paper. However, the relation to the proof of our Theorem~\ref{entropy2} is rather close so that it seems worthwile to include it here. In particular, both proofs are based on Laptev's result \cite{laptev}.

\begin{proof}[ of Theorem~\ref{thm:l2smallballs}] First note that the operator $u u^*: L_2[0,1] \to L_2[0,1]$\,, where $u^*$ denotes the adjoint operator of $u$, is  given by
$$
(u u^* g)(t) = \Gamma(2\beta+1)\,t^\alpha \int_0^1 \frac{  x^\alpha}{(t+x)^{2\beta+1}}\, g(x)\, \dd x\,.
$$
The exact asymptotic behavior of its singular numbers $s_n(u u^*)$\, was found by Laptev \cite{laptev} (based on Widom \cite{widom}). He showed that
\begin{equation} \label{eqn:dualoperatorev}
-\log s_n(u u^*) \sim 2 \pi (\alpha-\beta)^{1/2} n^{1/2}\,.
\end{equation}
Similar operators were studied in \cite{widom}, \cite{rk}, and \cite{bl}. Since $s_n(uu^*)=s_n(u)^2$ , this implies
\begin{equation} \label{eqn:atgn}
-\log s_n(u) \sim \pi (\alpha-\beta)^{1/2} n^{1/2}.
\end{equation} 
It is well-known that, due to the Karhunen-Lo\`eve expansion, it suffices to study the small deviation behavior of the weighted sum
$$
-\log \pr{ \sum_{n=1}^\infty s_n(u) \xi_n^2 \leq \eps^2 },
$$
which can be treated with Theorem~2 in \cite{nazarov} (cf.\ Proposition~4.1 in \cite{kn10}) or Theorem~4.1 in \cite{br}.  
\end{proof}

\medskip

As explained above, we give an alternative proof in order to provide a fully analytic line of arguments for Theorem~\ref{entropy2}.

\begin{proof}[ of Theorem~\ref{entropy2}] 
Recall that we determined the behavior of the singular numbers $s_n(u)$ in (\ref{eqn:atgn}). Since we are in the Hilbert space setting, the singular numbers and the entropy numbers of $u$ coincide with those of the diagonal operator $D_\sigma$ in $\ell_2$, $(x_n)\mapsto (\sigma_n x_n)$, with weight sequence $\sigma_n=s_n(u)$. A result of Gordon, K\"{o}nig, and Sch\"{u}tt (Proposition~1.7 in \cite{gks}) says  that
$$
\sup_{n\geq 1} \left(2^{-k/n} ( \sigma_1\cdots \sigma_n)^{1/n} \right)\leq e_{k+1}(D_\sigma) \leq 6\sup_{n\geq 1} \left( 2^{-k/n} ( \sigma_1\cdots \sigma_n)^{1/n}\right).
$$
This implies that
$$
\log e_{k}(u) \sim \sup_{n\geq 1} \left(-\frac{k}{n} \log 2 + \frac{1}{n} \sum_{j=1}^n \log s_j(u)\right).
$$
Using the asymptotics of $(s_n(u))$ from (\ref{eqn:atgn}), we get
$$
-\log e_{k}(u) \sim \inf_{n\geq 1} \left(\frac{k}{n} \log 2 + \frac{\pi (\alpha-\beta)^{1/2} }{n} \sum_{j=1}^n j^{1/2}\right).
$$
Evaluating this expression we obtain the assertion.
\end{proof}

\subsection{Proofs of Theorems~\ref{entropyinfty} and~\ref{general}}

\begin{proof}[ of Theorem~\ref{general}] Clearly, since $C_\lambda[0,1]$ is compactly embedded in $L_2[0,1]$, the operator $u:E\to L_2[0,1]$ is compact, whence its entropy numbers tend to zero.

For $0<\delta\le 1$ and $t\in[0,1]$ we consider the interval $I_\delta(t):=[0,1]\cap[t-\delta, t+\delta]$\, and define for $g\in L_2[0,1]$,
the local averaging operator
$$
(P_\delta g)(t):=\frac{1}{|I_\delta(t)|}\int_{I_\delta(t)}g(x)\,\dd x\,.
$$
One can easily verify that the averaging operators $P_\delta$ map $L_2[0,1]$ in $C[0,1]$.

{\it Step 1:} We need the following simple norm estimates:
\begin{align}
&\quad \Vert P_\delta : L_2[0,1] \to C[0,1] \Vert \le \delta^{-1/2}; \label{eqn:wasnaa}\\
&\quad \Vert {\rm id} - P_\delta : C_\lambda[0,1] \to C[0,1] \Vert \le \delta^{\lambda}. \label{eqn:wasnbb}
\end{align}
For all $t\in[0,1]$ and $g\in L_2[0,1]$ we have, by the Cauchy-Schwarz inequality,
$$
|(P_\delta g)(t)|\le \frac{1}{|I_\delta(t)|} \int_{I_\delta(t)} |g(x)|\,\dd x \le \frac{1}{|I_\delta(t)|}\, |I_\delta(t)|^{1/2}\Vert g\Vert_{L_2[0,1]}\le \delta^{-1/2}\Vert g\Vert_{L_2[0,1]}\,,
$$
which implies (\ref{eqn:wasnaa}).

Now let $t\in[0,1]$ and $g\in C_\lambda [0,1]$ with $\Vert g \Vert_{C_\lambda [0,1]}\le 1$. Then
$$
|g(t)-P_\delta g(t)|=\frac{1}{|I_\delta(t)|}\left|\int_{I_\delta(t)} (g(t)-g(x))\,\dd x \right|
\le \frac{1}{|I_\delta(t)|}\int_{I_\delta(t)} |g(t)-g(x)|\,dx\le \delta^\lambda\,,
$$
since $|g(t)-g(x)|\le |t-x|^\lambda \Vert g \Vert_{C_\lambda [0,1]}\le \delta^\lambda$ for every $x\in I_\delta(t)$\, and \, $\frac{1}{|I_\delta(t)|}\int_{I_\delta(t)} \,\dd x=1$.
To finish the proof of (\ref{eqn:wasnbb}), we take the supremum over all $t\in[0,1]$ and $g\in C_\lambda[0,1]$ with $\Vert g \Vert_{C_\lambda[0,1]}\le 1$\,.

{\it Step 2:} Using elementary properties of entropy numbers (cf.\ \cite{Carl}) and the above norm estimates (\ref{eqn:wasnaa}) and (\ref{eqn:wasnbb}), 
we obtain for all $k\in\N$ and $0<\delta \le 1$
\begin{align*}
e_k(u:E\to C[0,1])\le & \,e_k(P_\delta u:E\to C[0,1]) +\Vert u-P_\delta u:E\to C[0,1]\,\Vert\\
\le & \,e_k(u:E\to L_2[0,1])\cdot\Vert P_\delta :L_2[0,1]\to C[0,1]\,\Vert\\
& +\Vert u:E\to C_\lambda[0,1]\Vert\cdot\Vert {\rm id}-P_\delta : C_\lambda[0,1]\to C[0,1] \,\Vert\\
\le & \,e_k(u:E\to L_2[0,1])\cdot \delta^{-1/2} + \Vert u:E\to C_\lambda[0,1]\Vert\cdot\delta^\lambda\,.
\end{align*}
Finally we choose $\delta$ such that
$
e_k(u:E\to L_2[0,1])\cdot \delta^{-1/2} = \Vert u:E\to C_\lambda[0,1]\Vert\cdot\delta^\lambda\,,
$
i.e.\
$$
\delta =\left(\frac{e_k(u:E\to L_2[0,1])}{\Vert u:E\to C_\lambda[0,1]\Vert}\right)^\frac{1}{\lambda +1/2}\,.
$$ 
Clearly $0<\delta \le 1$, and consequently this choice of $\delta$ gives the desired inequality
$$
e_k(u:E \to C[0,1])\le 2\,\Vert u:E \to C_\lambda [0,1]\,\Vert^{\frac{1}{2\lambda +1}}\,e_k(u:E \to L_2[0,1])^{\frac{2\lambda}{2\lambda + 1}}\,.
$$
\end{proof}

\begin{proof}[ of the upper bound for the entropy numbers in Theorem~\ref{entropyinfty}]

{\it Step 1:} First we show that $u$ is a bounded operator from $L_2[0,\infty)$ into the H\"older space $C_\lambda[0,1]$, where
$$
\lambda:=\min(\alpha-\beta-1/2\,,\,1)\,.
$$ 

Let $f\in L_2[0,\infty)$, let $0\leq s<t \leq 1$, and set $h:=t-s$. We consider
$$
|(uf)(t)-(uf)(s)| \leq \int_0^\infty | t^\alpha x^\beta e^{- xt} - s^\alpha x^\beta e^{- xs}| \cdot |f(x)| \, \dd x.
$$
Using the Cauchy-Schwartz inequality, this can be estimated by $A^{1/2} \norm{f}_{L_2[0,\infty)}$, with
$$
A:= \int_0^\infty | t^\alpha x^\beta e^{- xt} - s^\alpha x^\beta e^{- xs}|^2 \dd x.
$$
We have to show that $A\leq C h^{2\lambda}$ for some constant $C$ independent of $t,s$, since then one can take the supremum over all $0\leq s<t \leq 1$ and finally over all $f$. Define
$$
g(y):=\int_0^\infty x^{2\beta} e^{-2yx} \dd x = y^{-2\beta-1} 2^{-2\beta-1}\Gamma(2\beta+1)
$$
and note that
\begin{align*}
A\,&=\int_0^\infty x^{2\beta} (t^{2\alpha} e^{- 2xt} -2 (ts)^\alpha e^{-x(t+s)} + s^{2\alpha} e^{- 2xs} ) \dd x\\
&= t^{2\alpha} g(t)  -2 (ts)^\alpha g((t+s)/2)+ s^{2\alpha} g(s)\,.
\end{align*}
Therefore, using the notation $\gamma:=\alpha-\beta-1/2$, we have 
\begin{align}
\frac{2^{2\beta+1}A}{\Gamma(2\beta+1)} &= t^{2\alpha-2\beta-1}  -2
(ts)^\alpha \left(\tfrac{t+s}{2}\right)^{-2\beta-1}+ s^{2\alpha-2\beta-1}
\label{eqn:htfs1}\\
&=\left[t^{2\gamma} -2 \left(\frac{t+s}{2}\right)^{2\gamma} + s^{2\gamma}\right] + 
2\left(\frac{t+s}{2}\right)^{-2\beta-1}\left[ \left(\frac{t+s}{2}\right)^{2\alpha}-
(ts)^\alpha\notag \right].
\end{align}

\noindent
{\it Case 1, $h\ge s$:} In this case we have $t=s+h\le 2s$ and $s\le h$, and therefore
(\ref{eqn:htfs1}) implies
$$
\frac{2^{2\beta+1}A}{\Gamma(2\beta+1)} \leq t^{2\gamma} + s^{2\gamma} \leq
(2h)^{2\gamma}+h^{2\gamma} = C h^{2\gamma}.
$$

\noindent {\it Case 2, $h\leq s$:} 
We estimate the two terms in the second line of (\ref{eqn:htfs1}) separately. For the first term 
we use the Taylor expansion of the function $f(x)=x^{2\gamma}$ at $x_0=(t+s)/2$ and obtain, with some $\xi\in(s,t)$, 
\begin{align}
t^{2\gamma} -2 \left(\frac{t+s}{2}\right)^{2\gamma} +
s^{2\gamma}&=f(x_0+h/2)-2f(x_0)+f(x_0-h/2)\label{1term}\\
&= \frac{h^2}{4} f''(\xi)=h^2\, \gamma(\gamma-1/2)\,\xi^{2\gamma-2}\le C_1h^2s^{2\gamma-2\notag }\,,
\end{align}
where we used $s\le \xi\le t=s+h\le 2s$ in the last inequality.

For the second term in (\ref{eqn:htfs1}) we apply the mean value theorem. Set $a:=ts$ and $b:=\left(\frac{t+s}{2}\right)^{2}$, and note that $a\leq b$ and $b-a=\frac{1}{4}(t^2+2ts+s^2-4ts)=(t-s)^2/4=h^2/4$. We obtain, with some $\eta\in (a,b)$, 
\begin{equation}\label{2term}
\left(\frac{t+s}{2}\right)^{2\alpha} - (ts)^\alpha
=(b-a)\cdot\alpha\, \eta^{\alpha-1}
\le C_2 h^2 s^{2(\alpha-1)}\,.
\end{equation}
Here we used that $s+t=2s+h\le 3s$, and therefore $s^2\le ts\le \eta\le (t+s)^2/4\le 9s^2/4$.

Combining now (\ref{eqn:htfs1}) with (\ref{1term}) and (\ref{2term}) we get
$$
\frac{2^{2\beta+1} A}{\Gamma(2\beta+1)}\leq C_1h^2s^{2\gamma-2}+2C_2s^{-2\beta-1}h^2s^{2\alpha-2}=Ch^2s^{2\gamma-2}\,.
$$
If $\gamma\ge 1$, then $h^2s^{2\gamma-2}\le h^2$; and if $\gamma <1$, we have
$h^2s^{2\gamma-2}=(h/s)^{2-2\gamma}h^{2\gamma}\le h^{2\gamma}$, since we are in the case $h\le s$. 
This proves, for all $0<h\le 1$, the desired estimate
$$
\frac{2^{2\beta+1} A}{\Gamma(2\beta+1)}\leq C h^{2 \min(\gamma, 1)}=Ch^{2\lambda}\,.
$$

{\it Step 2:} Let $\eps>0$. Then, by Theorem~\ref{entropy2}, for large enough $k$,
$$
e_k(u: L_2[0,\infty) \to L_2[0,1]) \leq e^{-(d_{\alpha,\beta}-\eps) k^{1/3}}.
$$
Now Theorem~\ref{general} implies that
$$
e_k(u: L_2[0,\infty) \to C[0,1]) \lesssim \exp\left(- \frac{2\lambda}{2\lambda +1}(d_{\alpha,\beta}-\eps) k^{1/3}\right)\,,
$$
and consequently
$$
-\log e_k(u: L_2[0,\infty) \to C[0,1]) \gtrsim \frac{2\lambda}{2\lambda +1}(d_{\alpha,\beta}-\eps) k^{1/3}\,.
$$
In other words, we have
$$
\liminf_{k\to\infty}\,\frac{-\log e_k(u: L_2[0,\infty) \to C[0,1]) }{k^{1/3}}\ge \frac{2\lambda}{2\lambda +1}(d_{\alpha,\beta}-\eps)\,.
$$
Letting $\eps \to 0$, we arrive at the desired upper bound for the entropy numbers.
\end{proof}

\begin{proof}[ of the lower bound for the entropy numbers in Theorem~\ref{entropyinfty}] 
Obviously, $e_k(u : L_2[0,\infty) \to C[0,1])\geq e_k(u : L_2[0,\infty) \to L_2[0,1])$. However, we can even gain a bit concerning the constant. 
For this purpose, let us stress the dependence on $\alpha$ and $\beta$ in the definition (\ref{uft}) by denoting the operator $u_{\alpha,\beta}$. 
Further, for some fixed $\eps>0$, we let $v : C[0,1]\to L_2[0,1]$ denote the multiplication operator $(vf)(t)=t^{-1/2+\eps} f(t)$. 
Note that $v : C[0,1] \to L_2[0,1]$ is bounded. Then one can observe that $v u_{\alpha,\beta} = u_{\alpha-1/2+\eps}$. Therefore,
$$
e_k(u_{\alpha-1/2+\eps} : L_2[0,\infty) \to L_2[0,1]) \leq e_k(u_{\alpha,\beta} : L_2[0,\infty) \to C[0,1])\cdot \norm{v : C[0,1] \to L_2[0,1]}.
$$
Using the $L_2$ estimate from Theorem~\ref{entropy2} for the left-hand side, this shows
$$
-\log e_k(u_{\alpha,\beta} : L_2[0,\infty) \to C[0,1]) \lesssim d_{\alpha-1/2+\eps,\beta} k^{1/3},
$$
or in other words,
$$
\limsup_{k\to\infty} \frac{-\log e_k(u_{\alpha,\beta} : L_2[0,\infty)\to C[0,1])}{k^{1/3}} \leq d_{\alpha-1/2+\eps,\beta}.
$$
This holds for all $\eps>0$. Letting $\eps$ tend to zero yields the lower bound for the entropy numbers, since the constant $d_{\alpha,\beta}$ (defined in (\ref{eqn:l2const})) is continuous in the parameters.
\end{proof}

\section{Relation to the entropy of function classes}\label{sec:functionS}
In this section, we relate the small deviation problem for $X_{\alpha,\beta}$ under the sup-norm to another small deviation problem, which in turn is related to a metric entropy problem of a certain function class. The function class related to the canonical case $\alpha=1$, $\beta=0$ is studied in Theorem~1.2 in \cite{GLW}. The present proof not only generalizes the case $\alpha=1$, $\beta=0$ but also gives a simplified proof for that specific case.

Let us define the process
$$
S(t):=t^{\alpha'}\int_0^1 x^\beta e^{-xt}\dd B(x),\qquad t\geq 1,
$$
where $\alpha':=2\beta+1-\alpha$ and the natural restrictions are $\alpha>\beta+1/2>0$. Note that exactly under these restrictions, $S$ is bounded on $[1,\infty]$.
Our main theorem concerning $S$ is as follows.
\begin{thm} \label{thm:S} Let $\alpha>\beta+1/2>0$. Then
\begin{equation} \label{eqn:thmS,s}
-\log\P \left( \sup_{t\geq 1} |S(t)| \leq \e \right) \asymp |\log \e|^3.
\end{equation}
\end{thm}

Using the technique in  \cite{GLW}, one finds that the associated convex hull for the process $S(t)$, $t \ge 1$, is the function class ${\cal F}$ consisting of all the functions $f$ on $[0,1]$ corresponding to the kernel $K(t,x)=t^{\alpha'}x^\beta e^{-tx}$. More precisely,  
${\cal F}$ can be expressed as
$$
{\cal F}:=\left\{f : [0,1] \to \R ~\left|~ f(x)=x^\beta \int_1^\infty t^{\alpha'} e^{-tx}\mu( {\rm d} t) \, : \,  \|\mu\|_{TV}\le 1 \right.\right\}.
$$
Under the $L_2[0, 1]$ norm $\|f\|_2=(\int_0^1 f^2(x) \dd x)^{1/2}$, the class ${\cal F}$ is compact and its metric entropy is denoted by $\log N(\e, {\cal F}, \|\cdot \|_2)$
where $N(\e, {\cal F}, \|\cdot \|_2)$ is the minimum number of $\e$-radii balls in the norm $\|\cdot \|_2$ to cover the class ${\cal F}$.
Thus, as discussed in detail in \cite{GLW}, via the connection between the small deviation probability and the metric entropy, 
we obtain the following statement for the function class ${\cal F}$ associated with $S$:

\begin{cor} For the class ${\cal F}$ defined above, and $\alpha'=2\beta+1-\alpha$ with $\alpha > \beta +1/2>0$,
\begin{equation} \label{eqn:cor-N}
\log N(\e, {\cal F}, \|\cdot \|_2) \asymp |\log \e|^3.
\end{equation}
\end{cor}

The original proof (as conducted in \cite{GLW} for the case $\alpha=1$, $\beta=0$) of the lower bounds for the estimates of the probability in (\ref{eqn:thmS,s}) and (\ref{eqn:supnormcase}) follows from covering estimates for the upper bound in  (\ref{eqn:cor-N}) which is lengthy and unpleasant.
The current approach for this part, which is turned around, is based on the simple and soft arguments summarized in Theorem~\ref{entropyinfty}.

However, the argument used for this -- that the upper bound of the metric entropy implies the lower bound of the small deviation probability, as discussed in Section~\ref{sec:entropymethod} -- is as in many other instances we know before: The key point is that it seems easier to find an upper bound of the metric entropy via analytic tools than a lower bound of the small deviation probability via probabilistic tools, even though they are equivalent. It would be interesting to find a probabilistic proof for the probability lower bound in (\ref{eqn:thmS,s}) or (\ref{eqn:supnormcase}) for all parameters in the range $\alpha >\beta+1/2>0$.

\begin{proof}[ of Theorem~\ref{thm:S}]
The class of processes that we consider satisfies the following time inversion property: $(X_{\alpha,\beta}(1/.))$ has the same law as $X_{\alpha',\beta}(.)$, which can be easily checked from the covariances.

For simplicity, set $\rho:=\sqrt{ \Gamma(2\beta+1) 2^{-(2\beta+1)} } $. Theorem~\ref{thm:linftysmallballs} yields that
$$
-|\log \e/\rho|^3\asymp \log\P\left( \sup_{t\leq 1} |\rho X_{\alpha,\beta}(t)|\leq \e \right) = \log\P\left( \sup_{t\geq 1} |\rho X_{\alpha',\beta}(t)|\leq \e \right).
$$
Clearly, by Anderson's inequality and the integral representation (\ref{def:x}), the last expression is smaller than
$$
\log\P\left( \sup_{t\geq 1} |S(t)|\leq \e \right),
$$
which already shows the lower bound of the small deviation probability of $S$ in (\ref{eqn:thmS,s}).

To see the opposite bound, note that
$$
\rho X_{\alpha',\beta}(t)= t^{\alpha'}\int_0^{2/\eps} x^\beta e^{-xt}\dd
B(x)+ t^{\alpha'}\int_{2/\eps}^\infty x^\beta e^{-xt}\dd B(x)=:V(t)+U(t)
$$
and thus
\begin{equation} \label{eqn:leftpr}
e^ {-c |\log \e|^3} \geq \P\left( \sup_{t\geq 1} |\rho
X_{\alpha',\beta}(t)|\leq \e \right) \geq \P\left( \sup_{t\geq 1} |V(t)|\leq
\e/2\right) \cdot  \P\left( \sup_{t\geq 1} |U(t)|\leq \e/2\right).
\end{equation}
Since
\begin{align*}
\P\left( \sup_{t\geq 1} |V(t)|\leq \e/2\right)
&= \P\left(\sup_{t\ge 1}\left|\int_0^{2/\eps}t^{\alpha'} x^\beta e^{-tx}\dd
B(x)\right|\leq \eps/2\right) \\
&= \P\left(\sup_{t\ge 1}\left|(\eps/2)^{\alpha'}\int_0^{1} (2t/\e)^{\alpha'}
(2x/\e)^\beta e^{-2x t/\e} (\e/2)^{-1/2} \dd B(x)\right|\leq \eps/2\right)
\\
&= \P\left(\sup_{t\ge 2/\e}\left|\int_0^{1} t^{\alpha'} x^\beta e^{-x t} \dd
B(x)\right|\leq (\eps/2)^{1+1/2+\beta-\alpha'}\right) \\
& \ge \P(\sup_{t\ge 1}|S(t)|\le (\eps/2)^{1+1/2+\beta-\alpha'})
\end{align*}
and $1+1/2+\beta-\alpha'>0$, it is sufficient to show that the second term
in (\ref{eqn:leftpr}) is bounded from below by a constant. This can be seen
as follows: Note that the finite dimensional distributions of $U$ are the same as of the following process \beaa & &t^{\alpha'}
e^{-t/\eps}\int_{1/\eps}^\infty (x+1/\eps)^\beta e^{-tx}\dd B(x)\\
&=&t^{\alpha'} \eps^{-\beta} e^{-t/\eps}\int_{1/\eps}^\infty (\eps
x+1)^\beta e^{-tx}\dd B(x)\\
&=& (t/\eps)^{\alpha'} \eps^{{\alpha'}-\beta-1/2}
e^{-t/\eps}\int_1^\infty (u+1)^\beta e^{-tu/\eps}\dd B(u).\eeaa
Therefore, estimating $e^{t/\eps}\geq e^{1/\eps}$ in the first step, we obtain
\beaa \P(\sup_{t\ge 1}|U(t)|\leq\eps/2)&\ge&
\P\left(\sup_{t\ge 1}\left|(t/\eps)^{\alpha'} \int_1^\infty
(u+1)^\beta
e^{-tu/\eps}\dd
B(u)\right|\leq\eps^{1+\beta+1/2-{\alpha'}}e^{1/\eps}/2\right)\\
&=& \P\left(\sup_{s\ge 1/\eps}\left|s^{\alpha'} \int_1^\infty
(u+1)^\beta
e^{-su}\dd B(u)\right|\leq\eps^{1+\beta+1/2-{\alpha'}}e^{1/\eps}/2\right)\\
&\ge& \P\left(\sup_{s\ge 1}\left|s^{\alpha'} \int_1^\infty
(u+1)^\beta e^{-su}\dd B(u)\right|\leq 1\right),
\eeaa
for small $\eps$ because $\eps^{1+\beta+1/2-{\alpha'}}e^{1/\eps}\to\infty$
as $\eps\to 0^+$. Note that the Gaussian process 
$$Z(s)=s^{\alpha'} \int_1^\infty
(u+1)^\beta e^{-su}\dd B(u)$$
is sample bounded on $[1,\infty)$ under the assumption $\alpha-\beta>1/2$.
Indeed,
\begin{eqnarray*}
\E
|Z(t)-Z(s)|^2&=&\int_1^\infty(u+1)^{2\beta}(t^{\alpha'}e^{-tu}-s^{\alpha'}e^{-su})^2\dd u\\
&\le& C\int_0^\infty
u^{2\beta}(t^{\alpha'}e^{-tu}-s^{\alpha'}e^{-su})^2 \dd u\\&=&C'\E|X_{\alpha',\beta}(t)-X_{\alpha',\beta}(s)|^2\\
&=&C'\E|X_{\alpha,\beta}(1/t)-X_{\alpha,\beta}(1/s)|^2.
\end{eqnarray*}
Now, Theorem~\ref{thm:linftysmallballs} implies that when $\alpha-\beta>1/2$, $Z(t)$ is sample
bounded on $[1,\infty)$. Therefore, $$\P\left(\sup_{s\ge 1}\left|s^{\alpha'}
\int_1^\infty
(u+1)^\beta e^{-su}\dd B(u)\right|\leq 1\right)$$
is a positive constant, as required.
\end{proof}


\begin{thebibliography}{10}
\bibitem{ab} H.\ Alzer and C.\ Berg.
\newblock Some classes of completely monotonic functions.
\newblock {\em Ann.\ Acad.\ Sci.\ Fenn.\ Math.} 27 (2002), 445--460.

\bibitem{ad09} F.\ Aurzada and S.\ Dereich.
\newblock Universality of the asymptotics of the one-sided exit problem for integrated processes.
\newblock To appear in: {\em Annales de l'Institut Henri Poincar\'e (B) Probabilit\'es et Statistiques}, preprint available from: http://arxiv.org/abs/1008.0485


\bibitem{ailz}
F.\ Aurzada, I.A.\ Ibragimov, M.A.\ Lifshits and J.H.\ van Zanten.
\newblock Small deviations of smooth stationary {G}aussian processes.
\newblock {\em Teor.\ Veroyatn.\ Primen.} 53 (2008), 788--798 (Russian); translation in {\em Theory Probab.\ Appl.} 53 (2009), 697--707 (English).



\bibitem{bl}
E.\ Belinsky and W.\ Linde.
\newblock Compactness properties of certain integral operators related to fractional integration.
\newblock {\em Math.\ Z.} 252  (2006),  669--686.



\bibitem{bgl} R.\ Blei, F.\ Gao and W.V.\ Li.
\newblock Metric entropy of high dimensional distributions.
\newblock  {\em  Proc. Amer. Math. Soc.} 135  (2007), 4009--4018.


\bibitem{br}
A.A.\ Borovkov and R.S.\ Ruzankin.
\newblock On small deviations of series of weighted random variables.
\newblock {\em J.\ Theoret.\ Probab.}  21  (2008),  628--649.

\bibitem{Carl}
B.\ Carl and I.\ Stephani.
\newblock {\em Entropy, compactness and the approximation of operators},
  volume~98 of {\em Cambridge Tracts in Mathematics}.
\newblock Cambridge University Press, Cambridge, 1990.



\bibitem{g04} F.\ Gao.
\newblock Entropy of absolute convex hulls in Hilbert spaces.
\newblock {\em Bull.\ London Math.\ Soc.} 36 (2004), 460--468.

\bibitem{g08} F.\ Gao.
\newblock  Entropy estimate for k-monotone functions via small ball probability
of integrated Brownian motion.
\newblock {\em Electron.\ Commun.\ Probab.} 13 (2008), 121--130.

\bibitem{GLW} F.\ Gao, W.V.\ Li and J.\ Wellner.
\newblock How many Laplace transforms of probability measures are there?
\newblock {\em Proc. Amer. Math. Soc.} 138 (2010), 4331--4344.

\bibitem{gw} F.\ Gao and J.\ Wellner.
\newblock  On the rate of convergence of the maximum
likelihood estimator of a k-monotone density.
\newblock {\em Science in China, Series A: Mathematics} 52 (2009), 1525--1538.

\bibitem{gks}
Y.\ Gordon, H.\ K\"{o}nig and C.\ Sch\"{u}tt.
\newblock Geometric and probabilistic estimates for entropy and approximation numbers of operators.
\newblock {\em J.\ Approx.\ Theory} 49 (1987), 219--239.

\bibitem{kn10}
A.I.\ Karol' and A.I.\ Nazarov.
\newblock Small Ball Probabilities for Smooth Gaussian fields and Tensor Products of Compact Operators.
\newblock Preprint available from: http://arxiv.org/abs/1009.4412


\bibitem{rk}
H.\ K\"{o}nig and S.\ Richter.
\newblock Eigenvalues of integral operators defined by analytic kernels.
\newblock {\em Math.\ Nachr.} 119 (1984), 141--155.

\bibitem{MR1701596}
C.\ Krattenthaler.
\newblock Advanced determinant calculus.
The Andrews Festschrift (Maratea, 1998).
\newblock {\it S\'{e}m.\ Lothar.\ Combin.} 42 (1999), Art. B42q, 67 pp. (electronic).

\bibitem{kl}
J.\ Kuelbs and W.V.\ Li.
\newblock Metric entropy and the small ball problem for {G}aussian measures.
\newblock {\em J.\ Funct.\ Anal.} 116 (1993), 133--157.

\bibitem{ku} T.\ K\"{u}hn.
\newblock Covering numbers of Gaussian reproducing kernel Hilbert spaces.
\newblock {\em J.\ Complexity} 27 (2011) 489--499.


\bibitem{laptev}
A.A.\ Laptev.
\newblock Spectral asymptotic behavior of a class of integral operators.
\newblock {\em Mat.\ Zametki} 16  (1974), 741--750 (Russian); translation in: {\em Math.\ Notes} 16 (1974), 1038--1043 (English).

\bibitem{Li10} W.V. Li.  Small Value Probabilities: Techniques and Applications. Lecture notes.
\verb+http://www.math.udel.edu/~wli/svp.html+


\bibitem{ll}
W.V.\ Li and W.\ Linde.
\newblock Approximation, metric entropy and small ball estimates for {G}aussian
  measures.
\newblock {\em Ann.\ Probab.} 27 (1999), 1556--1578.


\bibitem{LS}
W.V.\ Li and Q.-M.\ Shao.
\newblock \emph{Gaussian processes: inequalities, small ball
  probabilities and applications}, Stochastic processes: theory and methods,
  Handbook of Statist., vol.~19, pp.~533--597, North-Holland, Amsterdam, 2001.

\bibitem{LS04} W.V.\ Li and Q.-M.\ Shao.
\newblock Lower tail probabilities for Gaussian processes.
\newblock \emph{Ann.\ Probab.} 32 (2004), 216--242.

\bibitem{LS05} W.V.\ Li and Q.-M.\ Shao. Recent developments on lower tail probabilities for Gaussian processes.
\newblock \emph{COSMOS, the Journal of the Singapore National Academy of Science} 1 (2005), 95--106. 


\bibitem{sdbib} M.A.\ Lifshits.
  \newblock \emph{Bibliography compilation on small deviation probabilities},
   \newblock available from:\\ \verb+http://www.proba.jussieu.fr/pageperso/smalldev/biblio.html+,
  2010.

\bibitem{Majumdar}
S.N.\ Majumdar.
\newblock Persistence in nonequilibrium systems.
\newblock {\em Curr.\ Sci.}, 77 (1999), 370--375.
\newblock {\tt http://arxiv.org/abs/cond-mat/9907407}.


\bibitem{nazarov}
A.\ Nazarov.
\newblock Log-level comparison principle for small ball probabilities.
\newblock {\em Stat.\ Probab.\ Letters} 79 (2009), 481--486.

\bibitem{pietsch} A.\ Pietsch.
\newblock {\em Eigenvalues and $s$-numbers}, volume 13 of {\it Cambridge Studies in Advanced Mathematics}.
\newblock Cambridge University Press, Cambridge, 1987.

\bibitem{vaartzanten}
A.W.\ van der Vaart and J.H.\ van Zanten.
\newblock Rates of contraction of posterior distributions based on Gaussian process priors.
\newblock {\em Ann.\ Statist.} 36 (2008), 1435--1463. 

\bibitem{widom}
H.\ Widom.
\newblock Asymptotic behavior of the eigenvalues of certain integral equations. II.
\newblock {\em  Arch.\ Rat.\ Mech.\ Anal.}  17  (1964), 215--229.

\end{thebibliography}
\end{document}